\documentclass[]{article}
\usepackage{amssymb}
\usepackage{amsmath}
\usepackage{amssymb}
\usepackage{graphicx,epsfig}

\usepackage{epsfig}
\newtheorem{thm}{Theorem}
\newtheorem{pro}[thm]{Proposition}
\newtheorem{lem}[thm]{Lemma}
\newtheorem{cl}[thm]{Claim}
\newtheorem{cj}[thm]{Conjecture}
\newtheorem{cor}[thm]{Corollary}

\newcommand {\cbdo}{\hfill$\Box$}

\begin{document}

\title{Packing large trees of consecutive orders}
\author{Andrzej \.{Z}ak\thanks{The author was partially supported
by the Polish Ministry of Science and Higher Education.}\\
\small{AGH University of
Science and Technology, Krak\'ow,
Poland}} \maketitle

\begin{abstract}
A conjecture by Bollob\'as  from 1995 (which is a weakenning 
of the famous Tree Packing Conjecture by Gy\'arf\'as from 1976)  
states that any set of $k$ trees $T_n,T_{n-1},\dots,T_{n-k+1}$, 
such that $T_{n-i}$ has $n-i$ vertices, pack into $K_n$, 
provided $n$ is sufficiently large. We confirm Bollob\'as conjecture 
for trees $T_n,T_{n-1},\dots,T_{n-k+1}$, such that  
$T_{n-i}$ has $k-1-i$ leaves or a pending path of order $k-1-i$. 
As a consequence we obtain that the conjecture is true for $k\leq 5$. 
\end{abstract}

%%%%%%%%%%%%%%%%%%%%%%%%%%%%%%%%%%%%%%%%%%%%%%%%
\section{Introduction}

A set of (simple) graphs $G_1,G_2,\dots,G_k$
are said to \emph{pack into
a complete graph} $K_n$ (in short pack) if $G_1,G_2,\dots,G_k$
can be found as pairwise edge-disjoint subgraphs in
$K_n$. Many classical  problems in Graph Theory can be stated
as packing problems. In particular, $H$ is a subgraph of $G$ 
if and only if $H$ and the complement of $G$ pack. 

A famous tree packing conjecture
(TPC) posed by Gy\'arf\'as \cite{Gya} states that 
any set of $n$ trees $T_n,T_{n-1},\dots,T_1$ such that $T_i$ has $i$ vertices 
pack into $K_n$. A number of partial results concerning the TPC are known. 
In particular Gy\'arf\'as and Lehel \cite{Gya} showed that the TPC is true if each tree is either a path or a star. 
An elegant proof of this result was given by Zaks and Liu 
\cite{ZL}. 
Recently, B\"{o}ttcher et al. \cite{BHPT} proved an asymptotic version of the TPC for
trees with bounded maximum degree (see also \cite{MRS} for generalizations on other families of graphs).
In \cite{Bol2} Bollob\'as  suggested the following weakening of TPC 
\begin{cj}
For  every $k\geq 1$ there  is  an
$n_0(k)$
such  that  if
$n> n_0(k)$,  then
any set of
$k$
trees
$T_n,T_{n-1},\dots,T_{n-k+1}$
such that
$T_{n-j}$
has
$n-j$ vertices pack into $K_n$.
\end{cj} 
Bourgeois, Hobbs and Kasiraj \cite{HBK} showed that any three trees $T_n$, $T_{n-1}$, $T_{n-2}$ pack into $K_n$. 
Recently, Balogh and Palmer \cite{BP} proved
that any set of
$k=\frac{1}{10}n^{1/4}$ trees $T_n,\dots, T_{n-k+1}$ 
such that no tree is a star and
$T_{n-j}$ 
has
$n-j$ vertices pack into $K_n$. 
In this paper we confirm the conjecture for new sets of trees. 

We say that 
a tree $T$ has a \emph{pending path} of order $t$ if there exists  $e\in E(T)$ 
such that one component of $T-e$ is a path $P$ of order $t$ and $d_{T}(v)\leq 2$ for every $v\in V(P)$.  
\begin{thm}\label{main}
Let $k$ be a positive integer and let 
$n_0(k)$ be a sufficiently large constant depending only on $k$. 
If $n> n_0(k)$, then any set of $k$ trees $T_n,T_{n-1},\dots, T_{n-k+1}$, such that $T_{n-j}$ 
has $n-j$ vertices, and $T_{n-j}$ has $k-1-j$ leaves or a pending path of order $k-1-j$, 
pack into $K_n$.
\end{thm}
As an immediate consequence we obtain the following corollary.
\begin{cor}
Let $k\leq 5$ be a positive integer and let 
$n_0(k)$ be a sufficiently large constant depending only on $k$. 
If $n> n_0(k)$, then any set of $k$ trees $T_n,T_{n-1},\dots, T_{n-k+1}$, such that $T_{n-j}$ 
has $n-j$ vertices pack into $K_n$.
\end{cor}
The proofs of preparatory Lemmas \ref{pierwszy} and \ref{drugi} are inspired by 
Alon and Yuster approach \cite{AY}, 
but are much more involved. 

In what follows we fix an integer $k\geq 1$ and assume that $n\geq n_0(k)$, where $n_0(k)$ is a 
sufficiently large constant depending only on $k$.
\section{Notation}
The notation is standard. In particular $d_G(v)$ (abbreviated to $d(v)$ if no 
confusion arises) denotes the degree of a vertex $v$ in $G$, $\delta(G)$ and $\Delta(G)$ denote the minimum and the maximum degree of $G$, respectively.
Furthermore, $N_G(v)$ denotes the set of neighbors of $v$ and, for a subset 
of vertices $W\subseteq V(G)$, 
\[N_G(W)=\bigcup_{w\in W}N_G(w)\setminus W\] and 
\[N_G[W]=N_G(W)\cup W.\] 
Let $G$ be a graph and $W$  any set with $|V(G)|\leq |W|$. Given an 
injection $f:V(G)\rightarrow W$, 
let $f(G)$ denote the graph defined as fallows 
\[f(G)=\left(W,\{f(u)f(v):uv\in E(G)\}\right).\]  
For two graphs $G$ and $H$ let $G\oplus H$ denote the graph 
defined by 
\[G\oplus H=\left(V(G)\cup V(H),E(G)\cup E(H)\right)\] 
(note that $V(G)$ and $V(H)$ do not need to be disjoint). 

A packing of $k$ graphs $G_1,\dots,G_k$ with $|V(G_j)|\leq n$, $j=1,\dots,k$,  into 
a complete graph $K_n$ is a set of $k$ injections $f_j: V(G_j)\rightarrow V(K_n)$, $j=1,\dots,k$ such that 
\begin{align*}
& \text{if } i\neq j \text{ then } E(f_i(G_i))\cap E(f_j(G_j))=\emptyset.
\end{align*}
For two graphs $G$ and $H$ with $|V(G)|\leq |V(H)|$, we somtimes use an alternative definition. 
Namely, we call an injection $f:V(G)\rightarrow V(H)$ a packing of $G$ and $H$, if 
$E(f(G))\cap E(H)=\emptyset$. 
%%%%%%%%%%%%%%%%%%%%%%%%%%%%%%%%%%%%%%%%%%%%
\section{Preliminaries}
We write Bin$(p,n)$ for the binomial distribution with
$n$ trials and success probability $p$. Let $X\in {\rm Bin}(n,p)$. 
We will use the following two versions of the Chernoff bound which follows 
from formulas (2.5) and (2.6) from \cite{JLR} by taking 
$t=2\mu-np$ and $t=np-\mu/2$, respectively. 

If $\mu \geq E[X]=np$ then 
\begin{align}\label{chernoff1}
Pr[X\geq 2\mu]\leq {\rm exp}(-\mu/3)
\end{align}
On the other hand, if $\mu \leq E[X]=np$ then 
\begin{align}\label{chernoff2}
Pr[X\leq \mu/2]\leq {\rm exp}(-\mu/8).
\end{align}

\begin{pro}\label{ity}
Let $G$ be a graph with $n$ vertices and at most $m$ edges. Let $V(G)=\{v_1,\dots,v_n\}$ with 
$d(v_1)\geq d(v_2)\geq \cdots\geq d(v_n)$. Then 
\[d(v_i)\leq \frac{2m}{i}.\]
\end{pro}
Proof. The proposition is true because 
\[2m\geq \sum_{j=1}^{n}d(v_j)\geq \sum_{j=1}^id(v_j)\geq id(v_i).\]\cbdo
 
The following technical lemma is the main tool in the proof. A version 
of it appeared in \cite{AY}.  
\begin{lem}\label{cidi}
Let $G$ be a graph with $n$ vertices and at most $m$ edges. Let $V(G)=\{v_1,\dots,v_n\}$ with  
$d(v_1)\geq d(v_2) \geq \cdots \geq d(v_n)$. 
Let $A_i$, $i=1,\dots,n$, be any subsets 
of $V(G)$ with the additional requirement that if $u\in A_i$ then $d(u)< a$. 
For $i=1,\dots,n$ let $B_i$ be a random subset of $A_i$ where each vertex of $A_i$ is independently 
selected to $B_i$ with probability $p<1/a$. 
Let 
\begin{align*}
&C_i=\left( \bigcup_{j=1}^{i-1}B_j\right)\cap N(v_i),\\
&D_i=B_i\setminus \left( \bigcup_{j=1}^{i-1}N[B_j]\right).
\end{align*}
\noindent
Then  

1. $Pr\left[|C_i|\geq 4mp\right]\leq \exp(-2mp/3)$ for $i=1,\dots,n$

2. $Pr\left[|D_i|\leq \frac{p|A_i|}{2e}\right]\leq \exp\left( \frac{-p|A_i|}{8e}\right)$ for $i=1,\dots,\lfloor 1/(ap)\rfloor$.
\end{lem}
Proof. 
Fix some vertex $v_i\in V(G)$. 

Consider the first part of the lemma. 
If $d(v_i)\leq 2mp$ then the probability is zero because $|C_i|\leq |N(v_i)|=d(v_i)$. 
So we may assume that $d(v_i)> 2mp$. 
For $u \in N(v_i)$ 
the probability that $u\in B_j$ 
is at most $p$ (it is either $p$ if $u\in A_j$ or 0 if $u\not\in A_j$.) 
Thus $Pr[u\in C_i]\leq(i-1)p$. By Proposition \ref{ity},
$i \leq 2m/d(v_i)$. Hence, 
\[ Pr[u\in C_i]\leq \frac{2mp}{d(v_i)}.\]
Observe that
$|C_i|$ 
is a sum of $d(v_i)$ independent indicator random variables each of which 
has success probability at most $\frac{2mp}{d(v_i)}$. 
Thus, the expectation of $|C_i|$ is at most $2mp$. 
Therefore, by (\ref{chernoff1}), the probability
of
$|C_i|$ being larger than $4mp$ satisfies 
\begin{align*}
Pr[|C_i|\geq 4mp]\leq {\rm exp}\left( -2mp/3\right).
\end{align*}

Consider now the second part of the lemma. Observe that for $u\in A_i$, the probability that
$u\in B_i$ is $p$. 
On the other hand, for any
$j$, the probability that $u\not\in N[B_j]$ is at least $1-ap$. 
Indeed, 
$u\in N[B_j]$ if and only if $u\in B_j$ or one of its 
neighbors belongs to $B_j$. Since $u\in A_i$, it has at 
most $a-1$ neighbors. Hence, the probability that $u \in N[B_j]$ 
is at most $ap$.
Therefore, as long as $i\leq 1/(ap)$, 
\[Pr[u\in D_i]\geq p(1-ap)^{i-1}\geq \frac{p}{e}.\]
Observe that $|D_i|$ 
is a sum of $|A_i|$ 
independent indicator random variables, each having success
probability at least $\frac{p}{e}$. Therefore the expectation of $|D_i|$
is at
least 
$ \frac{p|A_i|}{e}$. 
By (\ref{chernoff2}), 
the probability that $|D_i|$ 
falls below $\frac{p|A_i|}{2e}$ satisfies
\begin{align*}
Pr\left[|D_i|\leq \frac{p|A_i|}{2e}\right]\leq {\rm exp}\left( -\frac{p|A_i|}{8e}\right).
\end{align*}
\cbdo
%%%%%%%%%%%%%%%%%%%%%%%%%%%%%%%%%%%%%%%%%%%%%%%%%%%%%%%%%%%%%%%%%%%%%%%%%%%%%%%
\section{Packing trees with small maximum degree.} 
\begin{lem}\label{pierwszy}
Let $G$ be a graph of order $n$ with $|E(G)|\leq kn$ and $\Delta(G)<2n/3+o(n)$. 
Let $T$ be a tree with $|V(T)|\leq n$ and  $\Delta(T)< 60(2k+1)n^{3/4}$. 
Let $I\subset V(G)$ with $|I|\leq k$ and such that if $v \in I$ then 
$d_G(v)\leq 2k$. Furthermore, let $I'\subset V(T)$ with $|I'|=|I|$ and such that 
if $v' \in I'$ then  $d_{T}(v')\leq 2$. 
Suppose, there is a packing $h':I'\rightarrow I$ of $T[I']$ and $G[I]$. 
Then, there is a packing $f':V(T)\rightarrow V(G)$ of $T$ and $G$ such that  
\begin{enumerate}
\item $\Delta(f'(T)\oplus G)\leq 2n/3+o(n)$,
\item $f'(v')=h'(v')$ for every $v'\in I'$.
\end{enumerate}
\end{lem}  
Proof. 
Let $V(G)=\{v_1,\dots,v_n\}$ where $d_G(v_i)\geq d_G(v_{i+1})$. Let $G'$ be a forest 
that arises from $T$ by adding $n-|V(T)|$ isolated vertices.  
Let $V(G')= \{v'_1,\dots,v'_{n}\}$ where $d_{G'}(v'_i)\geq d_{G'}(v'_{i+1})$. 
For convenience, we will construct a packing $f:V(G)\rightarrow V(G')$ 
such that $f(h'(v'))=v'$ for every $v'\in V(T)$. Thus for $f'$ we may take 
$f^{-1}$ restricted to $V(T)$.
 
Let $A_i\subset V(G)\setminus N_G[v_i]$ 
with the assumption that if $u\in A_i$ then $d_G(u)< 26k$.

\begin{cl}\label{sinowe}
$|A_i|\geq \frac{n}{4}$
\end{cl}
Proof. By the assumption on $\Delta(G)$, each vertex of $G$ has 
at least 
$n/3-o(n)$ 
non-neighbors. Suppose that $\alpha$ vertices of $G$ have degree 
greater than or equal to $26k$. Thus
\[ 2kn\geq 2|E(G)|=\sum_{i=1}^{n}d(v_i)\geq \alpha \cdot 26k,\]
and so 
$\alpha\leq \frac{ n}{13}$. Therefore 
\[|A_i|\geq n/3-o(n)-n/13\geq n/4.\]\cbdo \\

For $i=1,\dots,n$ let $B_i$ be a random subset of $A_i$ where each vertex of $A_i$ is independently 
selected to $B_i$ with probability 
\begin{align}\label{pe}
p=\frac{n^{-3/4}}{540\cdot26k^2(2k+1)}
\end{align}
Let 
\begin{align*}
&C_i=\left( \bigcup_{j=1}^{i-1}B_j\right)\cap N_G(v_i),\\
&D_i=B_i\setminus \left( \bigcup_{j=1}^{i-1}N_G[B_j]\right).
\end{align*}

\begin{cl}\label{cidi3}
The following hold with positive probability:

1. $|C_i|\leq \frac{n^{1/4}}{240(2k+1)}$ for $i=1,\dots,n$

2. $|D_i|\geq k(2k+1)+3$ for $i=1,\dots,\left\lfloor 540k(2k+1)n^{3/4}\right\rfloor$.
\end{cl}
Proof. Recall that $|E(G)|\leq kn$. 
Thus, by Lemma \ref{cidi}, the probability that $|C_i|>\frac{n^{1/4}}{240(2k+1)}$ ($>4|E(G)|p$), 
is exponentially small in $n^{1/4}$. 
Hence, for sufficiently large $n$ 
\[Pr\left[ |C_i|>\frac{n^{1/4}}{240(2k+1)}\right]<\frac{1}{2n}.\] 
Therefore, by the union bound, the first statement holds with probability greater than 1/2. 
Furthermore, by Claim \ref{sinowe},  
\[k(2k+1)+3<\frac{p|A_i|}{2e}.\]
Hence, by Lemma \ref{cidi} (with $a=26k$), for 
each $i\leq \left\lfloor 540k(2k+1)n^{3/4}\right\rfloor$ the probability 
that $|D_i|<k(2k+1)+3$ is exponentially small in $n^{1/4}$, as well. 
Hence, for sufficiently large $n$ 
\[Pr\left[ |D_i|<k(2k+1)+3\right]<\frac{1}{2n}.\] 
Therefore, by the union bound, the second statement holds with probability greater than 1/2, and so 
both statements hold with positive probability. \cbdo

Therefore, we may fix sets $B_1,\dots,B_n$ satisfying all the conditions 
of Claim \ref{cidi3} with respect to the cardinalities of the sets $C_i$ and $D_i$. 
We construct a packing $f:V(G)\rightarrow V(G')$ in three stages. At each point of the construction, 
some vertices of $V(G)$ are \emph{matched} to some vertices of $V(G')$, while the other 
vertices of $V(G)$ and $V(G')$ are yet \emph{unmatched}. Initially, all vertices are unmatched. 
We always maintain the packing property, 
that is for any $u,v\in V(G)$ if $uv\in E(G)$ then $f(u)f(v)\not\in E(G')$. 
The additional requirement that $\Delta(f(G)\oplus G')\leq 2n/3 +o(n)$ is preserved due 
to the assumption on $\Delta(T)=\Delta(G')$.  

After a forced Stage 1, in Stage 2 we match certain number of vertices of $G$ 
that have the largest degrees. After this stage, by the assumption 
on $\Delta(G')$, neither $G$ nor $G'$ has unmatched vertices of 
high degree (vertices of high degree are the main obstacle in packing). This fact 
enables us to complete the packing in Stages 3 and 4. 

\emph{Stage 1}
In Stage 1 we set $f(h'(v'))=v'$ for each $v'\in I'$. Clearly, 
the packing property is preserved.

\emph{Stage 2} Let $x$ be the largest integer such that $d_G(v_x)\geq \frac{n^{1/4}}{270(2k+1)}$. 
Thus, by Proposition \ref{ity},  
\begin{align}\label{ix2}
x\leq 540k(2k+1){n}^{3/4}
\end{align}

This stage is done repeatedly for $i=1,\dots, x$ and throughout it we maintain the following two invariants 
\begin{enumerate}
\item
At iteration $i$ we match $v_i$ with some vertex $f(v_i)$ of $G'$ 
such that $d_{G'}\left( f(v_i)\right)\leq 3$. 
\item
Furthermore,  
we also make sure that all neighbors of $f(v_i)$ in $G'$ are matched to vertices of $\bigcup_{j=1}^iB_j \cup I$. 
\end{enumerate}
To see that this is possible, consider the i'th iteration of Stage 1 where $v_i$ is some yet 
unmatched vertex of $G$. Let $Q'$ be the set of all yet unmatched vertices 
of $G'$ having degree less than or equal to 3. Note that, by Proposition \ref{ity}, the number of 
vertices of degree less than or equal to 3 in $G'$ is at least $n/2$. 
Hence, 
\[|Q'|\geq n/2-4(i-1)-k\geq n/2-4x-k\geq n/3.\]
Let $X$ be the set of already matched neighbors of $v_i$ and let 
$Y'=\bigcup_{u\in X}N_{G'}(f(u))$. 
Thus, the valid choice for $f(v_i)$ would be a vertex of $Q'\setminus Y'$. To see that 
such a choice is possible, it is enough to show that $|Q'|>|Y'|$.  
Let $X=X_1\cup X_2 \cup X_3$ with $X_1\subseteq I$, $X_2\subseteq \{v_1,\dots,v_{i-1}\}$ and 
$X_3\subseteq B_1\cup \cdots \cup B_{i-1}$. 
Hence $|X_1|\leq k$, $|X_2|\leq x$ and $|X_3|=|C_i|\leq \frac{n^{1/4}}{240(2k+1)}$. Thus, by the first invariant of 
Stage 2, and by (\ref{ix2}), Claim \ref{cidi3} and the assumptions on $I'$, 
\begin{align*}
|Q'|-|Y'|&\geq n/3-2|X_1|-3|X_2|-\Delta(G')|X_3|\geq n/3-2k-3x-60(2k+1)n^{3/4}\frac{n^{1/4}}{240(2k+1)}>0.
\end{align*}
In order to maintain the second invariant it remains to match 
the yet unmatched neighbors of $f(v_i)$ with vertices from $B_i$. 
Let $R'$ be the set of neighbors of $f(v_i)$ in $G'$ that are still unmatched. 
Recall that $|R'|\leq 3$. We have to match vertices of $R'$ with some vertices of $B_i$. 
Since $D_i=B_i\setminus \left( \bigcup_{j=1}^{i-1}N_G[B_j]\right)$, a valid choice of such 
vertices is by taking an $|R'|$-subset of $D_i\setminus N_G[I]$. By Claim \ref{cidi3} and by (\ref{ix2}), 
$|D_i|\geq k(2k+1)+3$ for $i=1,\dots,x$. Furthermore, since each $v\in D_i$ 
satisfies $d_G(v)<26k\leq d_G(v_x)$, $D_i\cap \{v_1,\dots,v_{i-1}\}=\emptyset$. 
Thus, the vertices from $D_i\setminus N_G[I]$ are still unmatched. Since $|N_G[I]|\leq k(2k+1)$ (by the assumptions on $I$), 
$|D_i\setminus N_G[I]|\geq 3$. Therefore, such a choice is possible.  

\emph{Stage 3}
Let $M_2$ and $M'_2$ be the set of matched vertices of $G$ and $G'$ after Stage 2, respectively. Clearly
$|M_2|=|M'_2|\leq 4x+k< n/9$. 
Hence $G'- M'_2$ has an independent set $J'$ with $|J'|\geq 4n/9$. 
Let $K'=V(G')\setminus (M'_2\cup J')$. 
In Stage 3 we match vertices of $K'$ one by one, with arbitrary yet unmatched vertices of $G$.  
Suppose that 
$v'\in K'$ is still unmatched. 
Let $Q$ be the set of all yet unmatched vertices of $G$. 
Clearly, $|Q|\geq |J'|\geq 4n/9$. Let $X'$ be the set of already 
matched neighbors of $v'$. Hence,  $|X'|\leq \Delta(G')\leq 60(2k+1)n^{3/4}$.  
Let $Y=\bigcup_{u'\in X'}N_{G}(f^{-1}(u'))$. 
 Thus, the valid choice for 
$f^{-1}(v')$ would be a vertex of $Q\setminus Y$. 
By the second invariant of Stage 2, 
$X'\cap \{v_1,\dots,v_x\}=\emptyset$. 
Hence, by the definition of $x$, 
\[|Y|\leq |X'|\cdot\frac{{n}^{1/4}}{270(2k+1)}\leq 60(2k+1)n^{3/4}\cdot\frac{{n}^{1/4}}{270(2k+1)}\leq 2n/9.\]
Therefore, $|Q\setminus Y|>0$, and so  an appropriate choice for $f^{-1}(v')$ is possible. 

\emph{Stage 4} 
Let $M_3$ and $M'_3$ be the sets of matched vertices of $G$ and $G'$ after Stage 3, respectively. 
In order to complete a packing of $G$ and $G'$, it remains to match 
the vertices of $V(G)\setminus M_3$ with the vertices of $J'$. 
Consider a bipartite 
graph $B$ whose sides are $V(G)\setminus M_3$ and $J'$. 
For two vertices $u\in V(G)\setminus M_3$ and $v'\in J'$, we 
place an edge $uv'\in E(B)$ if and only if it is possible 
to match $u$ with $v'$ (by this we mean that mapping $u$ to $v'$ 
will not violate the packing property). Thus $u$ is not allowed to be matched 
to at most $d_G(u)\Delta(G')$ vertices of $J'$. Thus  
\[d_B(u)\geq |J'|-\frac{{n}^{1/4}}{270(2k+1)}60(2k+1){n}^{3/4}\geq |J'|/2.\] 
On the other hand, since there is no edge 
from $v'$ to $v_i$ with $i\leq x $ (by the second invariant of Stage 2), 
$v'$ is not allowed to be matched to at most $\Delta(G')\frac{{n}^{1/4}}{270(2k+1)}$ vertices 
of $V(G)\setminus M_3$. Hence, analogously    
\[d_B(v')\geq  |J'|/2.\]
Therefore, by Hall's Theorem there is a matching of $V(G)\setminus M_2$ in $B$, and so a packing of $G$ and $G'$. 
\cbdo

%%%%%%%%%%%%%%%%%%%%%%%%%%%%%%%%%%%%%%%%%%%%
\section{Packing trees with large maximum degree}
\begin{lem}\label{drugi}
Let $G$ be a graph of order $n$ with $|E(G)|\leq kn$, $\delta(G)=0$ and $\Delta(G)<2n/3+o(n)$. 
Let $T$ be a tree with $|V(T)|\leq n$ and $\Delta(T)\geq 60(2k+1)n^{3/4}$. 
Let $I\subset V(G)$ with $|I|\leq k$ and such that if $v \in I$ then 
$d_G(v)\leq 2k$. Furthermore, let $I'\subset V(T)$ with $|I'|=|I|$ and such that 
if $v' \in I'$ then  $d_{T}(v')\leq 2$. 
Suppose, there is a packing $h':I'\rightarrow I$ of $T[I']$ and $G[I]$. 
Then, there is a packing $f':V(T)\rightarrow V(G)$ of $T$ and $G$ such that  
\begin{enumerate}
\item $\Delta(f'(T)\oplus G)\leq 2n/3+o(n)$,
\item $f'(v')=h'(v')$ for every $v'\in I'$.
\end{enumerate}
\end{lem}  
Proof. 
In the proof we will follow the ideas from the previous 
section. However, the key difference is that now both 
$G$  and $G'$ may have vertices of high degrees. Because of this obstacle,  
a packing has two more stages at the beginning. After a preparatory Stage 1,  
in Stage 2 we match the vertices of $G$ that have high degrees 
with vertices of $G'$ that have small degrees. Then in Stage 3, 
we match the vertices of $G'$ having high degree. This stage is very similar 
to Stage 2 from the previous section, but with the change of the role of 
$G$ and $G'$. Finally, we complete the packing in Stages 4 and 5, 
which are analogous to Stages 3 and 4 from the previous section.  

Let $V(G)=V=\{v_1,\dots,v_n\}$ where $d_G(v_i)\geq d_G(v_{i+1})$. 
Let $G'$ be a forest 
that arises from $T$ by adding $n-|V(T)|$ isolated vertices.  
Let $V(G')=V'= \{v'_1,\dots,v'_{n}\}$ where $d_{G'}(v'_i)\geq d_{G'}(v'_{i+1})$. 
For convenience, we will construct a packing $f:V\rightarrow V'$ 
such that $f(h'(v'))=v'$ for every $v'\in V(T)$. Thus for $f'$ we may take 
$f^{-1}$ restricted to $V(T)$. 

Let $A_i\subset V(G)\setminus N_G[v_i]$ 
with the assumption that if $u\in A_i$ then $d_G(u)< 26k$. 
The sets $A_i$ are defined in the same way as in the previous section. 
Thus,
\begin{align}\label{sinowe2}
|A_i|\geq \frac{n}{4}.
\end{align}
Let 
\begin{align}\label{defku}
q=\frac{{n}^{1/4}}{59(2k+1)}.
\end{align}
Let $P'\subseteq N_{G'}(v'_1)$ be the set of neighbors of $v'_1$ such that 
each vertex in $P'$ has degree at most $q$ in $G'$, and every  neighbor different from $v'_1$ of every vertex 
from $P'$  has degree at most $q$ in $G'$. 
\begin{cl}
$|P'|> (2k+1){n}^{3/4}$.
\end{cl}
Proof. Note that every vertex  $v'\in N_{G'}(v'_1)\setminus P'$ has the property that 
$d_{G'}(v')>q$ or $v'$ has a neighbor $w'\neq v'_1$ such that $d_{G'}(w')> q$. 
Therefore, 
\[n=|V(G')|> (\Delta(G')-|P'|)q\geq (60(2k+1){n}^{3/4}-|P'|)\frac{{n}^{1/4}}{59(2k+1)},\]
and the statement follows.
 \cbdo
\\
We construct a packing $f:V(G)\rightarrow V(G')$ in five stages. At each point of the construction, 
some vertices of $V(G)$ are \emph{matched} to some vertices of $V(G')$, while the other 
vertices of $V(G)$ and $V(G')$ are yet unmatched. Initially, all vertices are unmatched. 
We always maintain the packing property, 
that is for any $u,v\in V(G)$ if $uv\in E(G)$ then $f(u)f(v)\not\in E(G')$. 
Furthermore, we preserve that $\Delta(f(G)\oplus G')\leq 2n/3 +o(n)$.

\emph{Stage 1}. 
In Stage 1 we set $f(h'(v'))=v'$ for each $v'\in I'$. Furthermore we match 
an isolated vertex of $G$ with $v'_1$, i.e. $f(v_n)=v'_1$.

\emph{Stage 2}. Let $z$ be the largest integer such that 
$d_G(v_z)\geq {{n}^{1/4}}$. Since $|E(G)|\leq kn$, by Proposition \ref{ity} 
\begin{align}\label{ix}
z\leq  2k{n}^{3/4}.
\end{align}
This stage is done repeatedly for $i=1,\dots, z$ and throughout it we maintain the following invariants:
\begin{enumerate}
\item
At iteration $i$ we match $v_i$ with some vertex $f(v_i)$ of $G'$ 
such that $f(v_i)\in P'\setminus N_{G'}[I']$. 
\item Furthermore,  
we also make sure that all neighbors of $f(v_i)$ in $G'$, except $v'_1$, are matched to vertices of $A_i\setminus N_G[I]$. 
\end{enumerate}
Note that because $G'$ is a forest and since $P'\subseteq N_{G'}(v'_1)$, there are no edges 
between $N_{G'}[f(v_i)]$ and $N_{G'}[f(v_j)]$ for $i\neq j$. What is more, each $N_{G'}(f(v_j))$ is 
an independent set in $G'$. Since there are no edges (in $G$) between $v_i$ and $A_i$, 
the only edges that may spoil the packing property have one endpoint in $I$ or $I'$. 
However, by the first invariant there are no edges between $I'$ and $\bigcup_{j=1}^i f(v_j)$, 
and, by the second invariant, there are no edges between $I$ and $\bigcup_{j=1}^i N_G(v_i)\setminus\{v_1,\dots,v_i\}$. 
Therefore, such a mapping, if possible, do maintain the packing property. 
What is more, by (\ref{defku}) and by the definition of $P'$, 
the vertices of $G$ having large degrees are matched with vertices 
of $T$ having small degrees. Subsequently, by the definition of $z$, 
the vertices of $T$ having large degrees will be matched with vertices 
of $G$ having small degrees. Hence, the additional requirement that 
$\Delta(f(G)\oplus G')\leq 2n/3 +o(n)$ is preserved.

To see that this mapping is indeed possible, consider the $i$'th 
iteration of Stage 2, where $v_i$ is a vertex of 
$G$ with $d_G(v_i)\geq {n}^{1/4}\geq 26k$. In particular 
$v_i \not\in \bigcup_{j=1}^{i-1}A_j\cup I$, so  
$v_i$ is yet unmatched. 
Note that
\[|P'\setminus N_{G'}[I']|\geq (2k+1)n^{3/4}-3k\geq z\] 
and before iteration $i$, the number of 
already matched vertices of $P'\setminus N_{G'}[I']$ was equal to $i-1<z$.
Thus, there is at least one unmatched vertex in $P'\setminus N_{G'}[I']$, say $u'$, 
and we may set $f(v_i)=u'$ which preserves the first invariant.  
  
Furthermore, 
before iteration $i$ the overall number of matched vertices is at most 
\begin{align}\label{wzor1}
 |I|+1+(i-1)q < k+1+zq\leq k+n/59.
 \end{align}
Let $R'=N_{G'}(f(v_i))\setminus \{v'_1\}$. Note that all vertices from $R'$ are still unmatched. 
Thus, in order to maintain the second invariant, 
it suffices to match vertices of 
$R'$ with some vertices of $A_i\setminus N_G[I]$. 
Observe that by the choice of $P'$,  $|R'|\leq q-1$. 
Let $Q$ be the set of yet unmatched vertices of $A_i\setminus N_G[I]$.   
By (\ref{sinowe2}), (\ref{wzor1}), and since $|N_G[I]|\leq k(2k+1)$,  
\[|Q|\geq n/4-k(2k+1)-(k+n/59) > q-1.\] 
Hence, this is possible.  

Before we describe Stage 3, we need some preparations. 
Let $M_2$ be the set of all vertices of $G$ that were matched in Stage 1 or 2. 
Similarly, let $M'_2$ be the set of all vertices of $G'$ that were matched in Stage 1 or 2. 
Recall that 
\begin{align}\label{after2}
|M_2|=|M'_2|\leq k+1+zq<k+n/59.
\end{align}
Let $H=G[V\setminus M_2]$ be a subgraph of $G$ induced by yet unmatched vertices. 
Similarly let $H'=G'[V'\setminus M'_2]$. 
Note that since $G'$ is acyclic and by the construction of Stages 1 and 2, 
\begin{align}\label{stopienwh'}
&d_{G'}(v')\leq d_{H'}(v')+k+1 \text{ for each } v'\in V'\setminus M'_2.
\end{align}
Let $V(H')=\{w'_1,\dots,w'_r\}$ with $d_{H'}(w'_1)\geq d_{H'}(w'_2)\geq \cdots \geq d_{H'}(w'_r)$. 
By (\ref{after2}),
\begin{align}\label{er}
r\geq n-(k+n/59)>3n/4.
\end{align}
Let $y$ be the largest integer such that $d_{H'}(w'_y)\geq 360\sqrt{n}$. 
Then, by Proposition \ref{ity}, 
\begin{align}\label{igrek}
y\leq \frac{2n}{360\sqrt{n}}=\frac{\sqrt{n}}{180}.
\end{align}

For each $i=1,\dots,r$ we define a set 
$A'_i\subseteq V(H')\setminus N_{H'}[w'_i]$ to be a largest independent set of vertices but 
with the additional requirement that each $w' \in A'_i$ has $d_{H'}(w')< 180$.

\begin{cl}\label{siprim}
$|A'_i|\geq n/10$, $i=1,\dots,r$.
\end{cl}
Proof. Note that each $w'_i$ has at least 
\[r-d_{H'}(w'_i)-1\geq r-d_{G'}(w'_i)-1\geq r-d_{G'}(v'_2)-k-1\geq r-\frac{n}{2}-k-1\geq \frac{3}{4}n-\frac{n}{2}-k-1 =\frac{n}{4}-k-1\]
 non-neighbors. 
Since $H'$ is a forest, the subgraph of $H'$ induced by all non-neighbors of $w'_i$ 
has an independent set of cardinality at least $\frac{n/4-k-1}{2}>n/9$. 
Let $\alpha$ be the number of vertices of $H'$ that have degree 
greater than or equal to $180$. Thus
\[ 2n>\sum_{j=1}^{r}d_{H'}(w'_j)\geq \alpha \cdot 180,\]
and so 
$\alpha\leq \frac{n}{90}$. Therefore 
\[|A'_i|\geq n/9-\frac{n}{90}= n/10.\]
\cbdo

For $i=1,\dots,r$ let $B'_i$ be a random subset of $A'_i$ where each vertex of $A'_i$ is independently 
selected to $B'_i$ with probability $1/\sqrt{n}$. 
Let 
\begin{align*}
&C'_i=\left( \bigcup_{j=1}^{i-1}B'_j\right)\cap N_{H'}(w'_i),\\
&D'_i=B'_i\setminus \left( \bigcup_{j=1}^{i-1}N_{H'}[B'_j]\right).
\end{align*}

\begin{cl}\label{cidi2}
The following hold with positive probability:

1. $|C'_i|\leq 4\sqrt{n}$ for $i=1,\dots,r$

2. $|D'_i|\geq \frac{\sqrt{n}}{20e}$ for $i=1,\dots,y$.
\end{cl}
Proof. Clearly, $|E(H')|< n$. By Lemma \ref{cidi} 
(with $m=n$, $p=1/\sqrt{n}$ and $A_i=A'_i$), the probability 
that $|C'_i|\geq 4\sqrt{n}$ is exponentially small in $\sqrt{n}$. Thus, for $n$ sufficiently large  
\begin{align*}
&Pr\left[|C'_i|\geq 4\sqrt{n}\right]< \frac{1}{2n}\leq \frac{1}{2r}.
\end{align*}
Furthermore,  by Claim (\ref{siprim}),  
\[ \frac{\sqrt{n}}{20e}\leq\frac{p|A'_i|}{2e}.\]
Hence, by the second part of Lemma \ref{cidi}  (with $a=180$ and the remaining 
parameters as before) the probability that $|D'_i|\leq \frac{\sqrt{n}}{20e}$ 
is exponentially small in $\sqrt{n}$ for 
$i=1,\dots,\left\lceil \sqrt{n}/180 \right\rceil$. 
Thus,  by (\ref{igrek}), 
for $i\leq y\leq \lfloor\sqrt{n}/180\rfloor$ we have 
\begin{align*}
Pr\left[|D'_i|\leq \frac{\sqrt{n}}{20e}\right]&< \frac{1}{2y}.
\end{align*}
Thus, by the union bound,  each part of the lemma holds with probability greater than $1/2$. 
Hence both hold with positive probability. \cbdo

Now we are in the position to describe the next stages of a packing. By Claim \ref{cidi2} we may fix independent 
sets $B'_1,\dots,B'_r$ satisfying all the conditions of Claim \ref{cidi2} with respect to the cardinalities of 
the sets $C'_i$ and $D'_i$. 
Let $W=\{v_1,\dots,v_z\}$. Recall that 
\begin{align}\label{deltawu}
\Delta(G-W)<{{n}^{1/4}}.
\end{align}

\emph{Stage 3} 
This stage is done repeatedly for $i=1,\dots,y$ and throughout it we maintain the following 
two invariants 
\begin{enumerate}
\item
At iteration $i$ we match $w'_i\in V(H')$ with some yet unmatched vertex $u=f^{-1}(w'_i)$ of $H$ such that 
$d_G(u)\leq 4k$. 
\item
Furthermore,  
we also make sure that all neighbors of $f^{-1}(w'_i)$ in $H$ are matched to vertices of $\bigcup_{j=1}^iB'_j$. 
\end{enumerate}
To see that this is possible, consider the i'th iteration of Stage 3. Recall that $d_{H'}(w'_i)\geq 360\sqrt{n}\geq 180$. 
Hence, $w'_i$ does not 
belong to any $B'_j$ and so it is still unmatched. Let $Q$ be the set of all yet unmatched vertices 
of $G$ having degree less than or equal to $4k$. Note that, by Proposition \ref{ity}, the number of 
vertices of degree less than or equal to $4k$ in $G$ is at least $n/2$. 
Hence, by (\ref{after2}) and (\ref{igrek})
\begin{align}\label{ku}
|Q|&\geq n/2-|M_2|-(4k+1)y\geq n/2-k-n/59-(4k+1)\sqrt{n}/180>n/4.
\end{align}
Let $X'$ be the set of already matched neighbors in $G'$ of $w'_i$ 
and let $Y=\bigcup_{x'\in X'}N_{G}(f^{-1}(x'))$. 
Thus, the valid choice for $f^{-1}(w'_i)$ would be a vertex of $Q\setminus Y$. 
We will show that $|Q\setminus Y|>0$. 
Let $X'=X'_1\cup X'_2\cup X'_3$ such that $X'_1\subset M'_2$, $X'_2\subset \{w'_1,\dots,w'_{i-1}\}$ and 
$X'_3\subset \bigcup_{j=1}^{i-1}B'_j$. 
By (\ref{stopienwh'}), 
$|X'_1|\leq k+1$. Moreover if $v'\in X'_1$ then, by the 
second invariant of Stage 2, 
$v'\in M'_2\setminus \{f(v_1),\dots,f(v_z)\}$. Hence,  either $f^{-1}(v')\in I$ or $f^{-1}(v')$ belongs 
to some set $A_j$, $j\in \{1,\dots,z\}$.  Therefore, $d_G\left(f^{-1}(v')\right)\leq 26k$. 
Furthermore, $|X'_2|\leq i-1$ and, by Claim \ref{cidi2}, 
$|X'_3|=|C'_i|\leq  4 \sqrt{n}$. 
Hence, by (\ref{deltawu}) and by the first invariant of Stage 3,  
\begin{align}
|Y|&\leq 26k|X'_1|+4k|X'_2|+|X'_3|\cdot{{n}^{1/4}}
<n/4
\end{align}
Therefore, by (\ref{ku}),  
$|Q\setminus Y| >0$. 

In order to maintain the second invariant we have to match 
yet unmatched neighbors of $f^{-1}(w'_i)$ with some vertices 
of $B'_i$. 
Let $R$ be the set of the neighbors of $f^{-1}(w'_i)$ in $G$ that are still unmatched. 
Recall that $|R|\leq 4k$. 
Since $D'_i=B'_i\setminus \left( \bigcup_{j=1}^{i-1}N_{H'}[B'_j]\right)$, 
a natural choice of such 
vertices is by taking an $|R|$-subset of $D'_i$. 
However, unlike in Stage 2 in the previous subsection, 
this subset cannot be chosen arbitrarily because of 
the existence of possible edges between vertices from 
$P'':=I'\cup N_{G'}(P')\setminus \{v'_1\}$ and $D'_i$. 
For this reason, we have to match the vertices from $R$ more carefully. 
We match them, one by one,  
with some vertices from $D'_i$ in the following way. 
Suppose that $v\in R$ is yet unmatched. 
Let $D'$ be the set of yet unmatched vertices of $D'_i$. 
Since each $w'\in D'_i$ satisfies 
$d_{H'}(w')<180\leq 360\sqrt{n}$, $D'_i\cap\{w'_1,\dots,w'_{i-1}\}=\emptyset$. 
Hence, 
\begin{align}\label{debis}
|D'|\geq |D'_i|-|R|\geq {\sqrt{n}}/(20e)-4k.
\end{align}
Let $X_2$ be 
the set of all already matched neighbors of $v$ such that 
$f(X_2)\subseteq P''$. Let $Y'_2=\bigcup_{u\in X_2}N_{G'}(f(u))$. 
Thus, the valid choice for $f(v)$ would be a vertex from 
$D'\setminus Y'_2$.  Recall, that by the definition of $z$, 
$|X_2|\leq d_G(v)\leq n^{1/4}$. Furthermore, by the 
definition of $P'$ and $I'$, $|N_{G'}(f(u))|\leq q$. Thus, by (\ref{defku}) and (\ref{debis}),  
\begin{align*}
|D'\setminus Y'_2|>\sqrt{n}/(20e)-4k-|X_2|q\geq \sqrt{n}/(20e)-4k-\sqrt{n}/59>0.
\end{align*}
Thus, an appropriate choice for $f(v)$ is possible. 

\emph{Stage 4}
Let $M_3$ be the set of matched vertices of $G$ after Stage 3. 
Similarly, let $M'_3$ be the set of matched vertices of $G'$ after Stage 3. 
Note that, by (\ref{igrek}) and (\ref{after2}), 
\begin{align}\label{after3}
|M_3|=|M'_3|\leq |M_2|+(4k+1)y\leq k+n/59+(4k+1)\sqrt{n}/180<n/4
\end{align}
By (\ref{stopienwh'}), 
\begin{align}\label{deltawu'}
\Delta(G'-M'_3)\leq \Delta(H'-M'_3)+k+1\leq 360\sqrt{n}+k+1.
\end{align}
Furthermore, $|V(G')\setminus M'_3|> n-n/4=3n/4$. Thus $G'-M'_3$ has an independent set $J'$ with $|J'|> 3n/8$. 
Let $K'=V(G')\setminus (J'\cup M'_3)$. 
In Stage 4 we match vertices from $K'$ one by one, with arbitrary  yet unmatched vertices of $G$. 
Suppose that $v'\in K'$ is still unmatched. 
Let $Q$ be the set of all yet unmatched vertices of $G$. Clearly, $|Q|\geq |J'|\geq 3n/8$. 
Let $X'$ be the set of already matched neighbors of $v'$. 
By (\ref{deltawu'}), $|X'|\leq 360\sqrt{n}+k+1$. Let 
$Y=\bigcup_{x'\in X'}N_{G}(f^{-1}(x'))$. 
Thus, the valid choice for $f^{-1}(v')$ would be a vertex of $Q\setminus Y$.
By the second invariant of Stage 2, $X'\cap \{v_1,\dots,v_x\}=\emptyset$. 
Hence, by (\ref{deltawu}), 
\[|Y|\leq |X'|\cdot n^{1/4}<<3n/8-1.\]
Hence
\[|Q\setminus Y|\geq 1,\] 
and so the choice for $f^{-1}(v')$ is possible. 

\emph{Stage 5} 
Let $M_4$ and $M'_4$ be the sets of matched vertices of $G$ and $G'$, respectively, after Stage 4. 
In order to complete a packing of $G$ and $G'$ it remains to match 
the yet unmatched vertices of $G$ with vertices of $J'$. 
Consider a bipartite 
graph $B$ whose sides are $J:=V(G)\setminus M_4$ and $J'$. 
For two vertices $u\in J$ and $v'\in J'$, 
we place an edge $uv'\in E(B)$ 
if and only if it is possible to match 
$u$ with $v'$ (by this we mean that mapping $u$ to $v'$ will not violate the packing property). 
Recall that, by (\ref{deltawu}), $d_G(u)\leq {n}^{1/4}$. Moreover, by  the second invariant 
of Stage 3, $f(N_G(u))\subset V(G')\setminus \{w'_1,\dots,w'_y\}$. Thus, by 
the definition of $y$ and by (\ref{stopienwh'}), $u$ is not allowed 
to be matched to at most $n^{1/4}\left(360\sqrt{n}+k+1\right)$ vertices of $J'$. 
Therefore,  
\begin{align*}
d_B(u)\geq |J'|-{n}^{1/4}\left(360\sqrt{n}+k+1\right) 
> |J'|/2.
\end{align*}
Similarly, $d_{G'}(v')\leq 360\sqrt{n}+k+1$. 
Moreover, $f^{-1}(N_{G'}[v'])\subset V(G)\setminus W$. 
Thus, by (\ref{deltawu}), 
\[d_B(v')\geq |J'|-{n}^{1/4}\left(360\sqrt{n}+k+1\right)>|J'|/2.\] 
Therefore, by Hall's Theorem there is a perfect matching in $B$, and so a packing of $G$ and $G'$.
\cbdo 
%%%%%%%%%%%%%%%%%%%%%%%%%%%%%%%%%%%%%%%%%%%%%%%%%%%%%%%%%%%%%%%%%%%%%%%%%%%%%%%%%%%%%%%%
\section{Proof of Theorem \ref{main}}
Recall the theorem of Gy\'arf\'as and Lehel \cite{Gya}. 
\begin{thm}\label{pathstar}
Let $T_1,\dots,T_q$ be trees of orders $1,\dots,q$, respectively.  
If each $T_j$ is either a path or a star,  then there exists a packing of $T_j$, 
$j=1,\dots,q$, into $K_q$. 
\end{thm}

We will also need the following theorem proved by Brandt \cite{Br1}. 
\begin{thm}\label{brandt}
For every $0<\alpha<1/2$, there exists $n_0=n_0(\alpha)$ such that if $n> n_0$, 
$|E(G_1)|\leq \alpha n$ and 
$|E(G_2)|\leq \frac{1}{3\sqrt{\alpha}}n^{3/2}$, then $G_1$ and $G_2$ pack. 
\end{thm}

\noindent
{\bf Proof of Theorem \ref{main}.} 
We say that $T_{n-j}$, $j=0,\dots,k-1$, is of type I if $\Delta(T_{n-j})<60(2k+1)n^{3/4}$, of type II 
if  $60(2k+1)n^{3/4}\leq \Delta(T_{n-j})<2n/3$, and of type III if $2n/3\leq \Delta(T_{n-j})$. 
By the assumption, for each $j$ there exists a set $A_j\subset V(T_{n-j})$ such that 
either $A_j$ consists of $k-j-1$ leaves or $A_j$ is the vertex set of a pending path of order $k-j-1$. 
What is more, if $T_{n-j}$ is of type III, then it has a set $A_j$ which consists  
of $k-j-1$ leaf-neighbors of the vertex 
of maximum degree. 
If $A_j$ is the vertex set of a pending path then let $I_j=A_j\cup \{l\}$ where $l\in V(T_{n-j})\setminus A_j$ 
is a leaf 
of $T_{n-j}$ (clearly, such $l$ does exist). Otherwise let $I_j=A_j\cup \{w_j\}$ 
where $d_{T_{n-j}}(w_j)=\Delta(T_{n-j})$. 
Let $T_{k-j}=T_{n-j}[I_j]$. Thus $T_{k-j}=P_{k-j-1}\cup K_1$ or 
$T_{k-j}$ is a subgraph of the star $K_{1,k-j-1}$. 
In the former case let $T^*_{k-j}$ be a path that arises from $T_{k-j}$ by adding 
the missing edge beetwen $l$ and a vertex of degree $2$ in $T_{n-j}$. 
In the latter case, let $T^*_{k-j}$ be the star with the center $w_j$ 
that arises from $T_{k-j}$ by adding 
the missing edges incident to $w_j$. Note that if 
$T_{n-j}$ is of type III, then $T_{k-j}$ \emph{is} the star $K_{1,k-j-1}$.   
In particular for $j\in\{k-1,k-2,k-3\}$, $T^*_{k-j}$ is included to paths if $T_{n-j}$ is 
of type I, and $T^*_{k-j}$ is included to stars if $T_{n-j}$ is 
of type II or III. 
Clearly
\begin{align}\label{stst*}
d_{T_{n-j}}(v)\leq d_{T^*_{k-j}}(v) \text{ for every } v\in I_j\setminus w_j.
\end{align}
Let $G_0$ be a graph with vertex set $V=\{v_1,\dots,v_n\}$ and without edges.  
Let $K=\{v_{n-k+1},\dots,v_n\}$. By Theorem \ref{pathstar}, there exists a packing 
$h_j:V(T^*_{k-j})\rightarrow K$ of $T^*_{k-j}$, $j=0,\dots,k-1$. 

Let $T_{n-s_i}$, $i=1,\dots,s$, are of type III. 
We say that an edge $uv\in E(T^*_{k-j})$ is \emph{redundant} in $T^*_{k-j}$ 
(with respect to $h_j$, $j=0,\dots,k-1$),  
if 
\[\left|\{h_j(u),h_j(v)\}\cap \bigcup_{i=1}^s h_{s_i}\left(w_{s_i}\right)\right|\leq 1.\]
Otherwise $uv$ is called \emph{essential}.

A more detailed inspection of the short proof of Theorem \ref{pathstar} given bu Liu and Zaks \cite{ZL} (see also \cite{W}, p. 67), 
shows that there exists a packing $h_j: V(T^*_{k-j})\rightarrow K$, $j=0,\dots,k-1$, such that 
\begin{align}
%& \text{if } i\neq j \text{ then } E(h_i(T^*_{k-i}))\cap E(h_j(T^*_{k-j}))=\emptyset, \label{hstar0}\\ 
& \text{if } i>j, \text{ and } T^*_{k-i} \text{ and } T^*_{k-j} \text{ are stars, then }  
h_i(V(T^*_{k-i}))\subseteq K\setminus\{h_j(w_j)\} \label{hstar},\\
& \text{each path has a redundant edge incident to its endvertex}.\label{hstar2}
\end{align}
In particular, if $T_{n-j}$ is of type II then all edges 
of $T^*_{k-j}$ are redundant. Indeed, (\ref{hstar}) implies that 
$h_j(w_j)\neq h_{s_i}(w_{s_i})$, $i=1,\dots,s$. On the other hand all the edges of $T^*_{k-j}$ are 
incident to $w_j$.  
Furthermore, 
by (\ref{hstar2}), we may assume that if $T_{n-j}$ is of type I then $E(T^*_{k-j})\setminus E(T_{k-j})$ is a redundant edge. 
To sum up, we have that 
\begin{align}\label{redundant}
&T_{k-j}\text{ contains all essential edges of } T^*_{k-j} \;\;\;\;\;\;\;\;\;\;\;\;\;j=0,\dots,k-1.
\end{align}

Let $p,r,s$ be the numbers of trees of type I, II, and III, respectively. 
Let $P_1,\dots,P_p$ with $|P_i|<|P_{i+1}|$, $i=1,\dots,p-1$, denote the trees of type I. Similarly, let 
$R_1,\dots,R_r$ with $|R_i|<|R_{i+1}|$, and $S_1,\dots,S_s$, with $|S_i|<|S_{i+1}|$, denote the trees of 
type II and III, respectively. 
Let $R_i=T_{n-r_i}$, $i=1,\dots,r$. We partition $K$ and each $I_j$ into two subsets: 
\begin{align*}
&Y=\bigcup_{i=1}^{r}h_{r_i}(w_{r_i}),\\
&X=K\setminus Y,\\
&Y_j=h_j^{-1}(Y),\\
&X_j=h_j^{-1}(X)=I_j\setminus Y_j.
\end{align*} 
We first pack $R_i$, $i=1,\dots,r$, in a special way. 
We construct injections $f_{r_i}:V(R_i)\rightarrow V$, $i=1,\dots,r$, having the following properties: 
\begin{align*}
&E(f_{r_{i}}(R_i))\cap E(G_{i-1}) =\emptyset \text{ with } G_i=f_{r_i}(R_i)\oplus G_{i-1},\\
&f_{r_{i}}(v)=h_{r_{i}}(v) \text{ for every } v\in X_{r_i},\\
&f^{-1}_{r_i} (X)=X_{r_i},\\
&\Delta(G_{i})\leq 2n/3+o(n).
\end{align*}
To see that this is possible, consider the $i$-th iteration of this constructions. Note 
that, by (\ref{hstar}), 
$d_{G_{i-1}}(h_{r_i}(w_{r_i}))=0$. 
Hence $\delta(G_{i-1})=0$. 
Let $G'_{i-1}=G_{i-1}[V\setminus X\cup h_{r_i}(X_{r_i})]$. 
Since $h_{r_i}(w_{r_i})\in Y\subset V\setminus X$, $\delta(G'_{i-1})=0$, as well. 
Let $I'=X_{r_i}$ and 
$I=h_{r_i}(I')\subseteq X$. Clearly, $d_{R_i}(u)\leq 1$ for each $u\in I'$. 
Furthermore, $d_{R_t}(h_{r_t}^{-1}(v))\leq 1$, $t=1,\dots,i-1$,  for every $v\in I$. 
Hence, $d_{G'_{i-1}}(v)\leq i-1\leq k$ for every $v\in I$. 
Therefore, by Lemma \ref{drugi} with $G=G'_{i-1}$, $T=R_i$ and $h'=h_{r_i}$, 
an appropriate $f_{r_i}$ does exist. In particular, the third property is preserved because  
\begin{align*}
f_{r_i}^{-1}(X)&=f_{r_i}^{-1}(h_{r_i}(X_{r_i})) \;\;\;\;\;\;\;\;\;\;\;\;\;\;\text{ by the definition of } G'_{i-1}\\
&=f_{r_i}^{-1}(f_{r_i}(X_{r_i}))=X_{r_i} \;\;\;\;\text{ by the second property. }
\end{align*}

Now we pack $P_i:=T_{n-p_i}$, $i=1,\dots,p$. 
We construct injections $f_{p_i}:V(P_i)\rightarrow V$, $i=1,\dots,p$, having the following properties: 
\begin{align*}
&E(f_{p_{i}}(P_i))\cap E(G_{r+i-1}) =\emptyset \text{ with } G_{r+i}=f_{p_i}(P_i)\oplus G_{r+i-1},\\
&f_{p_{i}}(v)=h_{p_{i}}(v) \text{ for every } v\in X_{p_i},\\
&f^{-1}_{p_i} (X)=X_{p_i},\\
&\Delta(G_{i})\leq 2n/3+o(n).
\end{align*}
To see that this is possible, consider the $i$-th iteration of this constructions. 
Let $I'=X_{p_i}$ and 
$I=h_{p_i}(I')$. Clearly, $d_{P_i}(v)\leq 2$ for each $v\in I'\subseteq I_{p_i}$. 
Furthermore, $d_{P_t}\left(h_{p_t}^{-1}(v)\right)\leq 2$, $t=1,\dots,i-1$,  for every $v\in I\subseteq X$. 
Hence, $d_{G_{i-1}}(v)\leq r+2(i-1)\leq 2k$ for every $v\in I$. 
Therefore, by Lemma \ref{pierwszy} with $G=G_{r+i-1}[V\setminus X \cup h_{p_i}(X_{p_i})]$, $T=P_i$ and $h'=h_{p_i}$, 
an appropriate $f_{p_i}$ does exist. 

Finally, we pack $S_i:=T_{n-s_i}$, $i=1,\dots,s$. 
We distinguish in  $X$ and each $X_j$ the following subsets: 
\begin{align*}
&Z=\bigcup_{i=1}^{s}h_{s_i}(w_{s_i}),\\
&Z_j=h_j^{-1}(Z).
\end{align*} 
By (\ref{redundant}),  
\begin{align}\label{zjot}
&\left| N_{T_{k-j}}(v)\cap Z_j\right|=\left| N_{T^*_{k-j}}(v)\cap Z_j\right|
\end{align}
for every $v\in Z_j$, $j=0,\dots,k-1$.

Let $F_0=G_{p+r}$. 
Thus, 
\begin{align*}
F_0&=\bigoplus_{j\in[0,k-1]\setminus \{s_1,\dots,s_s\}}f_{j}(T_{n-j}).
\end{align*}
Let
\begin{align*}
H^*_{i-1}&=\bigoplus_{j\in[0,k-1]\setminus \{s_i,\dots,s_s\}}h_{j}(T^*_{k+1-j})\;\;\;\;\;\;\;\;\;\;i=1,\dots, s \\
H^*_{s}&=\bigoplus_{j\in[0,k-1]} h_{j}(T^*_{k+1-j}).
\end{align*}
Let $v_i=h_{s_i}(w_{s_i})$. By (\ref{hstar}), $v_i\neq v_j$ for $i\neq j$. 
We construct injections $f_{s_i}:V(S_i)\rightarrow V$, $i=1,\dots,s$, having the following properties: 
\begin{align*}
&E(f_{s_{i}}(S_i))\cap E(F_{i-1}) =\emptyset \text{ with } F_{i}=f_{s_i}(S_i)\oplus F_{i-1},\\
&f_{s_i}(v)=h_{s_i}(v) \text{ for every } v\in Z_{s_i},\\
&f^{-1}_{s_i} (Z) = Z_{s_i}. 
\end{align*}
To see that this is possible, consider the $i$-th iteration of this construction. 
We set $f_{s_i}(v)=h_{s_i}(v)$ for every  $v\in Z_{s_i}$. 
Note that $h_{s_i}$ sends $w_{s_i}$ and its $\alpha:=|Z_{s_i}|-1$ neighbors 
to $Z$. Since $w_{s_i}\in Z_{s_i}$, $f_{s_i}$ sends $w_{s_i}$ and its $\alpha$ neighbors 
to $Z$, as well. 
On the other hand $d_{T^*_{k-s_i}}=k-1-s_i$. 
Thus, $h_{s_i}$ sends $k-1-s_i-\alpha$ neighbors of $w_{s_i}$ to $K\setminus Z$. 
Since  
$h_j:V(T^*_{k-j})\rightarrow K$, $j=0,\dots,k-1$, is a packing 
of $T^*_{k-j}$, 
\begin{align}\label{ksi}
k-s-\left|N_{H^*_{i-1}}(v_i)\setminus Z\right|\geq k-1-s_i-\alpha.
\end{align}
We will show that there is a set $K'\subset (V\setminus Z)\setminus N_{F_{i-1}}(v_i) $ 
such that $|K'|=n-1-s_i-\alpha$. 
By (\ref{stst*}), (\ref{hstar}) and by the construction of $f_j$
\begin{align*}
\left|N_{F_{i-1}}(v_{i})\right|=d_{F_{i-1}}(v_{i})&=\sum_{j\in[0,k-1]\setminus\{s_i,\dots,s_s\},v_i\in f_{j}(Z_{j})}d_{T_{n-j}}\left(f^{-1}_j(v_i)\right)\\
&\leq \sum_{j\in[0,k-1]\setminus\{s_i,\dots,s_s\},v_i\in h_{j}(Z_{j})}d_{T^*_{k-j}}\left(h^{-1}_j(v_i)\right)=d_{H_{i-1}}(v_{i})=\left|N_{H^*_{i-1}}(v_{i})\right|.
\end{align*}
Therefore, by (\ref{zjot},\ref{ksi}), 
\begin{align*}
\left|(V\setminus Z)\setminus N_{F_{i-1}}(v_i)\right|&\geq n-s-\left|N_{F_{i-1}}(v_i)\setminus Z\right|
=n-s-\left|N_{F_{i-1}}(v_i)\right|+\left|N_{F_{i-1}}(v_i)\cap Z\right|\\
&\geq n-s-\left|N_{H^*_{i-1}}(v_i)\right|+\left|N_{H^*_{i-1}}(v_i)\cap Z\right|=
n-s-\left|N_{H^*_{i-1}}(v_i)\setminus Z\right|\\
&\geq n-1-s_i-\alpha.
\end{align*}
Thus, an appropriate $K'$ does exist.

Let $S'_{i}=S_i-Z_{s_i}$ and let $G'=F_{i-1}[K']$. 
Thus, $n':=|V(S'_i)|=n-s_i-\alpha-1$. 
In order to complete the construction of $f_{s_i}$, it is sufficient to pack 
$S'_i$ and $G'$. 
By the definition of trees of type III, $d_{S_i}(w_{s_i})\geq 2n/3$. 
Thus, $|E(S'_i)|\leq n/3\leq \frac{2}{5}n'$. 
Moreover, 
\[|E(G')|\leq kn << \frac{1}{3\sqrt{2/5}}(n')^{3/2}.\]
Thus, by Theorem \ref{brandt} 
such a packing does exist.\cbdo
\section{Remarks}
In the previous version of this paper, we claimed that we proved 
Bollob\'as conjecture in full. Unfortunately, the proof contained a mistake. 
We sincerely appologize the Readers for this misinformation.  

%----------------------------------------------------------------------------

\end{document}